\documentclass[12pt]{article}

\usepackage{amssymb}
\usepackage{amsmath}
\numberwithin{equation}{section}

\newtheorem{theorem}{Theorem}

\newtheorem{lemma}{Lemma}

\newcommand{\ooo}{\overline}

\newcommand{\ppp}{\partial}
\newcommand{\www}{\widetilde}

\newcommand{\N}{\Bbb{N}}
\newcommand{\R}{\Bbb{R}}

\title
{Determination of source terms and coefficients dependent on
$n-1$ spatial variable of parabolic and Schr\"odinger equations}

\author{$^1$ Oleg Imanuvilov,  
$^2$ Masahiro Yamamoto
}
\date{}
\begin{document}
\maketitle
\pagestyle{myheadings}
\begin{abstract}
{\it
We establish the uniqueness in the determination of a source term or
a coefficient of the zeroth order term of a second-order parabolic equation.
Moreover we consider the  determination of a potential of the
Schr\"odinger equation.
For a parabolic equation, an unknown source term and coefficient depend
on the time and  $n-1$ spatial variables.
For the Schr\"odinger equation, we assume that
a potential depends on $n-1$ spatial variables. The data are taken on a part
of the boundary satisfying some geometric condition. }
\end{abstract}


\section{Inverse Problem.}
The paper is concerned with the following two inverse problems:
\begin{itemize}
\item
for a second-order
for a parabolic equation: determination of a source term on the right-hand
side and a coefficient of the zeroth order term.
\item
(ii) for the Schr\"odinger equation: determination of a potential.
\end{itemize}

The data are collected on a part of the boundary.
For the parabolic case, unknown functions depend on the time variable 
and $n-1$ spatial
variables. Provided that an observation subboundary satisfies
some non restrictive and
easily verifiable geometric  condition,
we will prove the uniqueness for the determination of a source term and
a zeroth-order coefficient of the parabolic equation.

Now we formulate our inverse problems.
Let $\Omega$ be a bounded domain in $\Bbb R^n$ with $C^2-$ boundary and
$\nu =(\nu_1,\dots,\nu_n)$ be the outward unit normal vector to
$\partial\Omega$. Let $\Gamma$ be an arbitrarily fixed subboundary
of $\partial\Omega$.
We set
$Q=(0,T)\times \Omega$, $\Sigma=(0,T)\times \partial\Omega$,
$x=(x_1,\dots, x_n), x'=(x_1,\dots,x_{n-1})$.

For our inverse problems, 
we introduce an elliptic operator and related conditions:
Let functions $a_{ij}$ be real-valued
for all $i,j \in \{1,\dots,n\}$ and
\begin{equation}\label{sok4}
a_{ij}=a_{j i} \quad \mbox{for all $i,j\in \{1,\dots,n\}$}
\end{equation}
and
\begin{equation} a_{ij},\partial_{x_n}a_{ij} \in C^{1}(\ooo Q), \quad b_{i},
\partial_{x_n}b_i\in C^{0}(\ooo Q),\quad c,
\partial_{x_n}c\in  C^{0}(\ooo Q) \quad \mbox{for all $i,j\in \{1,\dots,n\}$}
\end{equation}
and there exists a constant $\beta>0$ such that
\begin{equation}\label{sok5}
\sum_{i,j=1}^n a_{ij}(t,x)\eta_{i}\eta_j\ge \beta\vert \eta\vert^2
\quad \mbox{for all $(t, x,\eta)\in Q\times \Bbb R^n$}.
\end{equation}
Furthermore, assume that there exists a constant $\beta_1>0$ such that
\begin{equation}\label{sok6}
R, \partial_{x_n}R\in C^{0}(\ooo Q), \quad \vert R(t,x)\vert >\beta_1>0\quad
\mbox{on }\,\, \ooo Q.
\end{equation}

Let $u(t,x)$ solve the following boundary value problem
\begin{equation}\label{Ssok1}
P(t,x,D)u:=\partial_tu-\sum_{i,j=1}^n \partial_{x_i}(a_{ij}(t,x)
\partial_{x_j}u)
+\sum_{i=1}^n b_{i}(t,x)\partial_{x_i} u+c(t,x) u=R(t,x) f(t,x')\quad
\mbox{in}\quad Q,
\end{equation}
\begin{equation}\label{SSsok2}
u\vert_{\Sigma}=0.
\end{equation}

Henceforth we set
$\partial_{\nu}u(t,x):= \sum_{i,j=1}^n a_{ij}(t,x)\partial_iu(t,x)\nu_j(x)$
for $(t,x) \in (0,T)\times \partial\Omega$.

Here a function $R$ is given, while a source term $f(t,x')$ is unknown.
Let $\Gamma_1$ be an open subset of $\Gamma$. 

We consider\\
{\bf Inverse source problem:}
\\
{\it Suppose that $(\partial_{\nu} u\vert_{(0,T)\times \Gamma}, \,
f\vert_{(0,T)\times \Gamma_1})$ are
given. Determine a source term $f$ in $Q$.}
\\

Let $\mathcal O$ be the projection of $\Omega$ on the
hyperplane $\{x\in \Bbb R^n\vert x_n=0\}$:
$$
\mathcal O=\{ x'\in \Bbb R^{n-1}\vert\, \mbox{there exists 
$y \in \Bbb R^1$ such that $(x',y)\in \Omega$}\}.
$$

We assume:
\\
{\bf Condition 1.} \label{fantomas}
{\it  There exists a dense set $\mathcal D$ in $\mathcal O$
such that for any $y'\in \mathcal D$, there exist $K(y') \in \N$, 
$y_n=y_n(y')\in \R$ satisfying $(y',y_n)\in \Gamma_1$, 
and connected open sets $S_k(y',y_n)\subset \Gamma$ 
for $k\in \{1,\dots, K(y')\}$
such that $S_{k-1}(y',y_n) \cap S_k(y',y_n) \ne \emptyset$ for
$2\le k \le K(y')$, $S_1(y',y_n) \cap \Gamma_1 \ne \emptyset$ and
$(y',y_n)\in S_{K(y')}(y',y_n)$.
Moreover for each $k \in \{ 1, 2,..., K(y')\}$, there exist a 
function $g_k(x') \in C^3(\overline{ O_k(y')})$
and a bounded domain $O_k(y') \subset \R^{n-1}$ with smooth boundary 
such that the set $S_k(y',y_n)$ is the graph of the function $g_k$:
$$
S_k(y',y_n)=\{ (x',g_k(x'))\vert x'\in O_k(y')\}.
$$
}

Then we can state our first main result.
\begin{theorem}\label{lox}
Let $(u_1,f_1),\, (u_2,f_2)\in H^{1,2}(Q)\times L^2(Q)$ satisfy 
(\ref{Ssok1}) and (\ref{SSsok2}), and let 
$\partial_{x_n}u_1,\partial_{x_n}u_2\in H^{1,2}(Q)$ and let
(\ref{sok4})- (\ref{sok6}) and Condition 1 hold.
If 
\begin{equation}\label{gavno}
(\partial_{\nu} u_1\vert_{(0,T)\times \Gamma}, \,
f_1\vert_{(0,T)\times \Gamma_1})
= (\partial_{\nu} u_2\vert_{(0,T)\times \Gamma}, \, 
f_2\vert_{(0,T)\times \Gamma_1}),
\end{equation}
then $f_1=f_2$ in $Q$.
\end{theorem}

Furthermore, we consider also an inverse problem of determining
a coefficient of zeroth order term depending on the time and
$n-1$ spatial variables.
Let $u_1(t,x)$ and $u_2(t,x)$ satisfy second-order parabolic equations
with zeroth order terms $c_1(t,x')$ and $c_2(t,x')$ respectively:
\begin{equation}\label{PKSsok1}
\partial_tu_1-\sum_{i,j=1}^n \partial_{x_i}(a_{ij}(t,x)\partial_{x_j} u_1)
+\sum_{i=1}^n b_{i}(t,x)\partial_{x_i} u_1+c_1(t,x') u_1=0\quad
\mbox{in}\quad Q,
\end{equation}
\begin{equation}\label{MKSsok1}
\partial_tu_2-\sum_{i,j=1}^n \partial_{x_i}(a_{ij}(t,x)\partial_{x_j} u_2)
+\sum_{i=1}^n b_{i}(t,x)\partial_{x_i} u_2+c_2(t,x') u_2=0\quad \mbox{in}
\quad Q
\end{equation}
and
\begin{equation}\label{KKSSsok2}
u_1=u_2\quad \mbox{on}\quad \Sigma.
\end{equation}
We introduce the following set
$$
\tilde \Psi=\mbox{Int}\{ x\in \Gamma_1 \vert\,\, c_1(t,x')=c_2(t,x')\quad
\mbox{for all $t\in (0,T)$} \}.
$$

We introduce a slightly different condition on $(y', y_n)$ and
$S_1(y',y_n)$ from Condition 1:
\\
{\bf Condition 2.} \label{fantomas}
{\it  There exists a dense set $\mathcal D$ in $\mathcal O$
such that for any $y'\in \mathcal D$, there exist $K(y') \in \N$, 
$y_n=y_n(y')\in \R$ satisfying $(y',y_n)\in \Gamma$, 
and connected open sets $S_k(y',y_n)\subset \Gamma$ 
for $k\in \{1,\dots, K(y')\}$
such that $S_{k-1}(y',y_n) \cap S_k(y',y_n) \ne \emptyset$ for
$2\le k \le K(y')$, $S_1(y',y_n) \cap \tilde \Psi \ne \emptyset$ and
$(y',y_n)\in S_{K(y')}(y',y_n)$.
Moreover for each $k \in \{ 1, 2,..., K(y')\}$, there exist a 
function $g_k(x') \in C^3(\overline{ O_k(y')})$
and a bounded domain $O_k(y') \subset \R^{n-1}$ with smooth boundary 
such that the set $S_k(y',y_n)$ is the graph of the function $g_k$:
$$
S_k(y',y_n)=\{ (x',g_k(x'))\vert x'\in O_k(y')\}.
$$
}

We have
\begin{theorem}\label{Plox}
Let $u_1, u_2\in H^{1,2}(Q)\cap C^0(\ooo{Q})$ satisfy (\ref{PKSsok1}) -
(\ref{KKSSsok2}) and $\partial_{x_n}u_1,\partial_{x_n}u_2\in H^{1,2}(Q)
\cap C^0(\ooo Q)$.  We assume (\ref{sok4})- (\ref{sok6}) and Condition 2,
and the existence of a constant $\beta>0$ such that
\begin{equation}\label{gavno1}
\vert u_2(t,x)\vert>\beta>0\quad \mbox{in}\quad Q
\end{equation} and
\begin{equation}\label{gavno}
(\partial_{\nu} u_1\vert_{(0,T)\times \Gamma}, \, 
c_1\vert_{(0,T)\times \Gamma_1})
= (\partial_{\nu} u_2\vert_{(0,T)\times \Gamma}, \,
c_2\vert_{(0,T)\times \Gamma_1}).
\end{equation}
Then $c_1=c_2$ in $Q$.
\end{theorem}

Next we consider an inverse problem of the determination of a potential
depending on $n-1$ spatial variables for a Schr\"odinger equation.
More precisely, $u_1(t,x)$ and $u_2(t,x)$ are assumed to satisfy the
Schr\"odinger equations with potentials $c_1$ and $c_2$ respectively:
\begin{equation}\label{XPKSsok1}
\sqrt{-1}\partial_tu_1-\sum_{i,j=1}^n \partial_{x_i}(a_{ij}(x)
\partial_{x_j} u_1)
+ \sum_{i=1}^n b_{i}(x)\partial_{x_i} u_1+c_1(x') u_1=0\quad \mbox{in}\quad Q,
\end{equation}
\begin{equation}\label{XMKSsok1}
\sqrt{-1}\partial_tu_2-\sum_{i,j=1}^n \partial_{x_i}(a_{ij}(x)
\partial_{x_j} u_2)+\sum_{i=1}^n b_{i}(x)\partial_{x_i} u_2+c_2(x') u_2=0
\quad \mbox{in}\quad Q
\end{equation}
and
\begin{equation}\label{XKSSsok2}
u_1=u_2\quad \mbox{on}\quad \Sigma.
\end{equation}
Here we assume that $a_{ij}, b_i, c_1, c_2$ satisfy (\ref{sok4}) - 
(\ref{sok5}), and the coefficients $b_i$ and $c_k$ are assumed to be 
complex-valued, not necessarily real-valued.
\\

We have
\begin{theorem}\label{XPlox}
Let $u_1, u_2\in H^{1,2}(Q)\cap C^0(\ooo Q)$ satisfy
(\ref{XPKSsok1})-(\ref{XKSSsok2}), and $\partial_{x_n}u_1,\partial_{x_n}u_2
\in H^{1,2}(Q)\cap C^0(\ooo Q)$ and let (\ref{sok4})-(\ref{sok6}) and
 Condition 2 hold true.
Moreover, we assume that there exists a constant $\beta>0$ such that
\begin{equation}\label{Pgavno1}
\vert u_2(t,x)\vert>\beta>0\quad \mbox{on}\quad Q.
\end{equation}
If 
\begin{equation}\label{Pgavno}
(\partial_{\nu} u_1\vert_{(0,T)\times \Gamma}, \,
c_1\vert_{(0,T)\times \Gamma_1})
=(\partial_{\nu} u_2\vert_{(0,T)\times \Gamma}, \,
c_2\vert_{(0,T)\times \Gamma_1}).
\end{equation}
Then $c_1=c_2$ in $Q$.
\end{theorem}

The inverse problem of the determination of a coefficient independent on
one component of $(x,t)$ of parabolic
equation was studied in \cite{Bez1} in the case of $\Omega = \Bbb R^{n}$
under some restrictive geometry of $\Gamma$ and $\Gamma_1$.
Moreover it is assumed that an initial value is known.
The paper \cite{Bez2} discusses the case of $\Omega=\Bbb R^n_+$.
In both papers, the author derived integral equations with respect to
unknown coefficients.
The papers \cite{GK} and \cite{IKY} established the uniqueness and
a conditional stability estimate for the determination of a source term and
a coefficient in case where $\Omega$ is cylindrical and the principal term of
the parabolic equation is either $\Delta$ or the coefficients of the principal
part are $x_n$-independent.

The current paper generalizes the uniqueness result of \cite {IKY} in the
sense that we consider a general parabolic equation and a general
domain $\Omega.$

The problem of an inverse problem of determining a time independent potential
in the Schr\"odinger equation was studied in Bukhgeim \cite{BU}.
The observation data are given on a part of the boundary which
satisfies some non-trapping type of geometrical conditions, provided that
data of some trace at some moment of time are known.
We refer to recent works \cite{B-P}, \cite{BP1}, \cite{IY5}, \cite{IY6},
where observation is taken at final or initial time moment.
Unlike these existing works, our Theorem \ref{XPlox} requires neither
spatial data over domain at a time moment nor non-trapping type
geometrical conditions on the observation part of the boundary,
but a potential is assumed to depend only  on $n-1$ spatial variable.

The proof is based on the Carleman estimates with non-zero boundary data,
and we modify the methodology by a pioneering work Bukhgeim and Klibanov
\cite{BK}.

The paper is composed of three sections.  Sections 2 and 3 are devoted to 
the proofs of Theorems 1-2 and Theorem 3 respectively. 

\section{Proof of Theorem \ref{lox}.}

The proof of Theorem \ref{lox} is based on the following lemma.
Let  functions $f$ and $u$ satisfy (\ref{Ssok1}) and (\ref{SSsok2}). 
We introduce
$$
\mathcal G=\mbox{Int}\,\{x \in \partial\Omega \vert f(t,x)=u(t,x)
=\nabla u(t,x)=0 \quad \mbox{for all $t\in (0,T)$} \}.
$$

Assume

{\bf Condition 3.} \label{fantomas}
{\it 
For $y'\in O$, we assume that there exist a point $(y',y_n)\in \Gamma$ and 
an open connected set $S(y',y_n)\subset \Gamma$ such that 
$S(y',y_n) \cap \mathcal G $ is not empty, $(y',y_n)\in S(y',y_n)$,
and we can find a function $g:= g(x') \in C^3(\overline{O(y')})$
such that the set $S(y',y_n)$ is a graph of $g:$
$$
S(y',y_n) =\{ (x',g(x'))\vert x'\in O(y') \}.
$$
Here $O(y')$ is a bounded domain in $\Bbb R^{n-1}$  with smooth boundary.}
\\

Then
\begin{lemma}\label{lox1} 
Let $(u,f)\in H^{1,2}(Q)\times L^2(Q)$ satisfy 
$\partial_{x_n}u\in H^{1,2}(Q)$, (\ref{Ssok1}) and (\ref{SSsok2}). 
Suppose that (\ref{sok4})- (\ref{sok6}) and  Condition 3 
hold true and
$$
(\partial_{\nu} u\vert_{(0,T)\times \Gamma}, \,
f\vert_{(0,T)\times \Gamma_1})=(0, 0).
$$ 
Then there exists a neighborhood $U\subset \Bbb R^{n-1}$ of $y'$ such that 
$f=0$ in $(0,T)\times U.$
\end{lemma}

{\bf Proof.}
We make a change of the variables in a neighborhood  $\Pi$ of
the surface $S(y',y_n)$:
$$
\widetilde x=F(x)=(x', x_n-g(x')).
$$
An unknown function $f$ is not essentially affected by this change of 
the variables.

On the other hand  there exist a strictly positive $\epsilon$  such that
$$
\widetilde \Omega=O(y')\times(0,\epsilon)\subset F(\Pi\cap \Omega).
$$
The change of coordinates transforms the parabolic operator $P$ into 
a parabolic operator $\widetilde P$:
\begin{equation}\label{focus}
\widetilde P(t,\widetilde x,D)\widetilde u
:= \partial_t\widetilde u-\sum_{i,j=1}^n \partial_{\widetilde x_i}
(\widetilde a_{ij}(t,\widetilde x)\partial_{\widetilde x_j}\widetilde u)
+ \sum_{i=1}^n \widetilde b_{i}(t,\widetilde x)\partial_{\widetilde x_i}
\widetilde{u}
\end{equation}
$$ 
+ \widetilde c(t,x) \widetilde u
= \widetilde R(t,\widetilde x) f(t,\widetilde x')
\quad \mbox{in}\quad (0,T)\times F(\Pi\cap \Omega),
$$
and
\begin{equation}\label{Ssok2}
\partial_{\widetilde x_n}\widetilde u\vert_{(0,T)\times O(y')\times\{0\}}=
\widetilde u\vert_{(0,T)\times O(y')\times\{0\}}=0.
\end{equation}

By (\ref{sok4})-(\ref{sok6}), the coefficients $\widetilde a_{ij},
\widetilde b_j,\widetilde c$ satisfy
\begin{equation}\label{ZKsok4}
\widetilde a_{ij}=\widetilde a_{j i}\quad \mbox{for all
$i,j\in \{1,\dots,n\}$} 
\quad \mbox{and}\quad \mbox{for all $(t,\widetilde x)\in (0,T)\times 
F(\Pi\cap \Omega)$}
\end{equation}
and
\begin{eqnarray} 
\widetilde a_{ij},\, \partial_{\widetilde x_n}\widetilde a_{ij} 
\in C^{1}(\overline{(0,T)\times F(\Pi\cap \Omega)}), 
\quad \widetilde b_{i},\, \partial_{\widetilde x_n}\widetilde b_i
\in C^{0}(\overline{(0,T)\times F(\Pi\cap \Omega)}),\nonumber\cr\\
\quad \widetilde c,\, \partial_{\widetilde x_n} \widetilde c
\in  C^{0}(\overline{(0,T)\times F(\Pi\cap \Omega)})\quad 
\mbox{for all $i,j\in \{1,\dots,n\}$}
\end{eqnarray}
and there exist constants $\beta>0$ and $\widetilde{\beta_1}>0$ such that
\begin{equation}\label{ZKsok5}
\sum_{i,j=1}^n \widetilde a_{ij}(t,\widetilde x)\eta_{i}\eta_j
\ge \beta\vert \eta\vert^2 \quad \mbox{for all 
$(t, \widetilde x,\eta)\in (0,T)\times F(\Pi\cap \Omega)\times \Bbb R^n$}
\end{equation}
and 
\begin{equation}\label{Ksok6}
\widetilde R, \, \partial_{\widetilde x_n}\widetilde R
\in C^{2}(\overline{(0,T)\times F(\Pi\cap \Omega)}), \quad 
\vert \widetilde R(t,\widetilde x)\vert >\widetilde \beta_1>0\quad 
\mbox{on }\,\, \overline{(0,T)\times F(\Pi\cap \Omega)}.
\end{equation}
Let us introduce a weight function
\begin{equation}\label{sous}
\psi(\widetilde x)=\widetilde \psi(\widetilde x')-N\widetilde x_n^2+K.
\end{equation}
Here
$N$ are $K$ are large positive parameters satisfying 
\begin{equation}\label{kkm}
\inf_{\widetilde x'\in O(y')}\psi(\widetilde x',0)
> \sup_{\widetilde x'\in O(y')}\psi(\widetilde x',\epsilon)
\end{equation} 
and 
\begin{equation}\label{Puk1} 
\psi(\widetilde x)>0\quad \mbox{on}\,\,\overline{\widetilde  \Omega}.
\end{equation}
Let $\widetilde \psi\in C^2(\overline{O(y')})$ satisfy 
\begin{equation} \label{Puk}
\left\{ \begin{array}{rl}
& \nabla\widetilde\psi\ne 0\quad\mbox{on}\,\, \overline{O(y')},\quad  
\partial_{\nu'}\widetilde \psi<0\quad \mbox{on} \quad \partial O(y'), 
                                               \cr\\
& \widetilde \psi\vert_{\partial O(y')\setminus F(\mathcal G)}=0, 
\quad \widetilde \psi(x')>0\quad \quad\mbox{on}\quad \overline{O(y')},
\end{array}\right.
\end{equation}
Here $\nu'$ denotes the unit outward normal vector to $\partial O(y')
\subset \R^{n-1}$ and $\ppp_{\nu'}\www{\psi} = \nabla_{x'}\www{\psi}
\cdot \nu'$.

The existence of such a function $\www{\psi}$ is established in \cite{Ima1995}.

For some  $t_0\in (0,T)$, let
\begin{equation}\label{klown}  
f(t_0,y')\ne 0.
\end{equation}
In view of $\psi$, we introduce two functions $\alpha(t,\widetilde x)$ 
and $\varphi(t,\widetilde x)$:
\begin{equation}\label{DD1.2} 
\varphi(t, \widetilde x) = \frac{e^{\lambda \psi(\widetilde x)}}{\ell(t)},
\quad \alpha(t, \widetilde x) 
= \frac{e^{\lambda\psi(\widetilde x)} 
- e^{2\lambda ||\psi||_{C^0(\overline{\widetilde \Omega})}}}
{\ell(t)},
\end{equation}
where $\ell(T)\in C^\infty [0,T]$, $\ell(0)=\ell(T)=0$, 
$\ell'(0)\ne 0$, $\ell'(T)\ne 0$, and $\ell$ is strictly monotone increasing  
in $(0,t_0-\delta_0]$ and strictly monotone decreasing in $[t_0+\delta_0,T]$ 
and $\ell(t)=1$ on  $[t_0-\delta_0,t_0+\delta_0]$ with a small constant 
$\delta_0>0$.
Therefore
\begin{equation}\label{klown1}
1=\ell(t_0)>\ell(t) \quad \mbox{for all 
$t\in (0,T)\setminus [t_0-\delta_0,t_0+\delta_0]$}.
\end{equation}
We set
$$
\alpha_I(t)=\inf_{\widetilde x'\in O(y')}\alpha(t,\widetilde x',0)
\quad \mbox{and}\quad \alpha_M(t)=\sup_{\widetilde x'\in  O(y')}
\alpha(t,\widetilde x',\epsilon).
$$
By (\ref{kkm}), there exists a constant $\delta>0$ such that
\begin{equation}\label{khmer}
\alpha_I(t)-\alpha_M(t)\ge \frac{\delta}{\ell(t)},
\quad \forall t\in [0,T].
\end{equation}

Let us consider the boundary value problem
\begin{equation}\label{LD1.4}
 \widetilde P (t,\widetilde x,D)z =  g 
\quad \mbox{in}\quad \widetilde Q:=(0,T)\times O(y')\times (0,\epsilon),
\end{equation}
\begin{equation}\label{LD1.5}
z \big|_{\widetilde \Sigma_0} = 0,\quad z(0, \cdot) = z_0,
\end{equation}
where
$$
\widetilde \Sigma_0=(0,T)\times O(y')\times \{0\}.
$$

In order to formulate the next lemma, we introduce two operators
\begin{eqnarray}\label{sosna}
\widetilde P_1(t,\widetilde x,D,\tau)z
= \partial_t z
+ 2\tau\varphi\lambda \tilde a(t,\widetilde x,\nabla \psi,\nabla z)\nonumber\\
+ \tau\lambda^2\varphi \widetilde a(t,\widetilde x,\nabla\psi,\nabla\psi)z
+ \tau\varphi\lambda \sum_{i,j=1}^n \partial_{\widetilde x_i}\widetilde 
 a_{ij}\partial_{\widetilde x_j}\psi z,
\end{eqnarray}
$$
\widetilde P_2(t,\widetilde x,D,\tau)z
=- \sum_{i,j=1}^n \widetilde a_{ij} \partial^2_{\widetilde x_i\tilde x_j} z 
- \lambda^2 {\tau}^2 \varphi^2 \widetilde a(t, \widetilde x, \nabla\psi, 
\nabla\psi) z,
$$ 
where $\widetilde a(t,\widetilde x,\eta,\eta)=\sum_{i,j=1}^n \widetilde 
a_{ij}(t,\widetilde x)\eta_i\eta_j,\nabla=(\partial_{\widetilde x_1},\dots, 
\partial_{\widetilde x_n})$.
Let $\gamma=\sum_{j=1}^n\Vert \widetilde b_j\Vert_{L^\infty(\widetilde Q)}
+ \Vert \widetilde c\Vert_{L^\infty(\widetilde Q)}$.
\\

Then
\begin{lemma}\label{D1.2} {\bf (Carleman estimate).}
{\it Let (\ref{ZKsok4})-(\ref{ZKsok5}) hold true, the functions 
$\varphi, \alpha$ be defined by (\ref{DD1.2}), and let 
$\psi$ satisfy (\ref{Puk1}) and(\ref{Puk}). 
Let $z\in H^{1,2}(\widetilde Q)$ and $g\in L^2(\widetilde Q)$.
Then there exists a constant
$\widehat\lambda(\gamma) > 0$ satisfying:  
for an arbitrary $\lambda \ge \widehat{\lambda}$ there exists
${\tau}_0(\lambda,\gamma) > 0$ such that solutions of problem 
(\ref{LD1.4}) - (\ref{LD1.5}) satisfy the following inequality:
\begin{eqnarray}\label{DD1.6}
  \int_{\widetilde Q} ( {\tau}\lambda^2\varphi |\nabla z|^2 
+ {\tau}^3\lambda^4 \varphi^3 \vert z\vert^2)
e^{2{\tau}\alpha}   d\widetilde x\,dt
+ \sum_{k=1}^2\Vert \widetilde P_k(t,\widetilde x,D,\tau) (ze^{\tau\alpha})
\Vert^2_{L^2(\widetilde Q)}             \nonumber\\
\le C_1\left( \int_{\widetilde Q} |g|^2 e^{2{\tau}\alpha} d\widetilde xdt 
+ \int_{\widetilde\Sigma\setminus\widetilde \Sigma_0} 
(\tau \lambda\varphi \vert\nabla z\vert^2+\tau^3 \lambda^3\varphi^3 
\vert z\vert^2)e^{2{\tau}\alpha}  d\widetilde \Sigma\right)
\end{eqnarray}
for each ${\tau} \ge \tau_0(\lambda,\gamma)$.
Here the constant $C_1$ is independent of $\tau$ and $\lambda.$
}
\end{lemma}

We use Lemma \ref{D1.2} to complete the proof of Lemma \ref{lox1}.

Applying Carleman estimate (\ref{DD1.6}) to the system (\ref{focus}),
we obtain
\begin{eqnarray}\label{P1}
\int_{\widetilde Q} ( {\tau}\lambda^2\varphi |\nabla {\widetilde u}|^2 
+ {\tau}^3 \lambda^4\varphi^3 \vert {\widetilde u}\vert^2)
e^{2{\tau}\alpha}  d\widetilde x\,dt
+ \sum_{k=1}^2\Vert\widetilde P_k(t,\widetilde x,D,\tau) 
({\widetilde u}e^{\tau\alpha})\Vert^2_{L^2({\widetilde Q})} \nonumber\\
\le C_2\left( \int_{\widetilde Q} |f|^2 e^{2{\tau}\alpha}  d\widetilde x dt 
+ \int_{\widetilde\Sigma\setminus\widetilde \Sigma_0} (\tau \lambda\varphi 
\vert\nabla {\widetilde u}\vert^2 + \tau^3 \lambda^3\varphi^3 \vert 
{\widetilde u}\vert^2)e^{2{\tau}\alpha}  d\widetilde\Sigma\right)
\end{eqnarray}
for all $ \tau\ge \tau_0(\lambda)$ and $ \lambda\ge \lambda_0 .$
Next we differentiate  both sides of equation (\ref{focus}) 
by $\partial_{\widetilde x_n}:$

\begin{equation}\label{focus1}
\widetilde P(t,\widetilde x,D)\partial_{\widetilde x_n}\widetilde u
= [\partial_{\widetilde x_n},\widetilde P] \widetilde u
+ \partial_{\widetilde x_n}\widetilde R(t,\widetilde x) f(t,\widetilde x')
\quad \mbox{in}\quad \tilde Q,
\end{equation}
\begin{equation}\label{fSsok2}
\partial_{\widetilde x_n}\widetilde u\vert_{(0,T)\times O(y')\times\{0\}}=0.
\end{equation}
Applying (\ref{DD1.6}) to (\ref{focus1}) and (\ref{fSsok2}), we obtain
\begin{eqnarray}\label{P11}
\int_{\widetilde Q} ( {\tau}\lambda^2\varphi |\nabla \partial_{\widetilde x_n}
{\widetilde u}|^2 \nonumber\\
+ {\tau}^3 \lambda^4\varphi^3 \vert \partial_{\widetilde x_n}{\widetilde u}
\vert^2) e^{2{\tau}\alpha} d\widetilde x\,dt
+ \sum_{k=1}^2\Vert \widetilde P_k(t,\widetilde x,D,\tau)
 ((\partial_{\widetilde x_n}{\widetilde u})e^{\tau\alpha})\Vert^2
_{L^2(\widetilde Q)}                                 \nonumber\\
\le C_3\left( \int_{\widetilde Q} (|f|^2 e^{2{\tau}\alpha}   
+ \sum_{\vert\gamma\vert\le 2}\vert \partial^\gamma_{\widetilde x} 
{\widetilde u}\vert^2 e^{2\tau\alpha}) d\widetilde xdt\right.\nonumber\\
+ \left.\int_{\widetilde\Sigma\setminus \widetilde \Sigma_0} 
(\tau \lambda \varphi \vert\nabla\partial_{\widetilde x_n}{\widetilde u}
\vert^2
+ \tau^3 \lambda^3 \varphi^3 \vert\partial_{\widetilde x_n}{\widetilde u}
\vert^2) e^{2{\tau}\alpha}  d\widetilde \Sigma\right)
\end{eqnarray}
for all $ \tau\ge \tau_1(\lambda)$ and $ \lambda\ge \lambda_1.$
From (\ref{P1}) and (\ref{P11}), we have
\begin{eqnarray}\label{P1Z}
\int_{\widetilde Q} ( {\tau}\lambda^2\varphi |\nabla \partial_{\widetilde x_n}
{\widetilde u}|^2 
+ {\tau}^3\lambda^4 \varphi^3 \vert \partial_{\widetilde x_n}{\widetilde u}
\vert^2) e^{2{\tau}\alpha} d\widetilde x\,dt
+ \sum_{k=1}^2\Vert \widetilde  P_k(t,\widetilde x,D,\tau) 
((\partial_{\widetilde x_n}{\widetilde u})e^{\tau\alpha})\Vert^2
_{L^2({\widetilde Q})}                  \nonumber\\
\le C_4\left( \tau\int_{\widetilde Q}\varphi |f|^2 e^{2{\tau}\alpha} 
d\widetilde xdt 
+ \int_{\widetilde\Sigma\setminus\widetilde \Sigma_0}\sum_{\ell=0}^1 
(\tau \lambda \varphi \vert\nabla\partial^\ell_{\widetilde x_n}{\widetilde u}
\vert^2
+ \tau^3 \lambda^3 \varphi^3 \vert\partial^\ell_{\widetilde x_n}{\widetilde u}
\vert^2) e^{2{\tau}\alpha} d\widetilde\Sigma\right)
\end{eqnarray}
for all $\tau\ge \tau_2(\lambda)$ and $ \lambda\ge \lambda_2$.
By (\ref{sous}) and (\ref{DD1.2}), we obtain
\begin{align*}
& \int_{\widetilde Q}\varphi |f|^2 e^{2{\tau}\alpha (t,\widetilde x)}
  d\widetilde x dt\\
=& \int_{(0,T)\times O(y') } \frac{e^{\lambda\widetilde\psi(\tilde x')+\lambda K}}
{\ell(t)} |f|^2  e^{-\frac{2{\tau}}{\ell(t)}e^{2\lambda\Vert \psi\Vert
_{C^0(\ooo{\widetilde\Omega})}}} 
\int_0^\epsilon  e^{-\lambda N\widetilde x_n^2}e^{\frac{2\tau }{\ell(t)} 
e^{\lambda\widetilde \psi(\widetilde x')+\lambda K } 
e^{-\lambda N \widetilde x_n^2}}d\widetilde x_n d\widetilde x'dt
\end{align*}
$$
= \int_{(0,T)\times O(y') } \frac{e^{\lambda\widetilde\psi(\tilde x')+\lambda K}}
{\ell(t)}|f|^2 e^{2{\tau}\alpha (t,\widetilde x',0)}
$$
$$
\times \left( \int_0^{\epsilon\,\root\of{\lambda N}}  
e^{-\widetilde x_n^2}
\exp\left( \frac{2\tau }{\ell(t)} 
e^{\lambda\tilde \psi(\widetilde x')+\lambda K} 
(e^{-\tilde x_n^2}-1)\right ) d\widetilde x_n\right) 
d\www{x}' 
\frac{dt}{\sqrt{\lambda N}}.
$$
%
Function $p(x_n)=e^{-x_n^2}-1\le 0, p(0)=p'(0)=0$ and $p''(0)=-2.$ Using the Morse lemma (see e.g. \cite{Her}) one can find a coordinates $x_n=r(y_n)$ such that $(e^{-\tilde x_n^2}-1)\circ r(y_n)=-y_n^2.$
Making the corresponding change of variables, one can estimate the following integral as 

\begin{eqnarray}\label{okorok}
 \int_0^{\epsilon\,\root\of{\lambda N}}  
e^{-\widetilde x_n^2}
\exp\left( \frac{2\tau }{\ell(t)} 
e^{\lambda\tilde \psi(\widetilde x')+\lambda K} 
(e^{-x_n^2}-1)\right ) d\widetilde x_n\nonumber\\\le C_5 \int_{-\infty}^{\infty}  
e^{-\widetilde x_n^2}
\exp\left( -\frac{2\tau }{\ell(t)} 
e^{\lambda\tilde \psi(\widetilde x')+\lambda K} 
y_n^2\right ) dy_n \le \frac{C_6}{\root\of{\tau\frac{2\tau }{\ell(t)} 
e^{\lambda\tilde \psi(\widetilde x')+\lambda K}} }.
\end{eqnarray}


Using (\ref{okorok}), we see that there exists a constant $C_5>0$ such that
\begin{equation}\label{okorok1}
\int_{\widetilde Q}\tau \varphi |f|^2 e^{2{\tau}\alpha (t,\widetilde x)}  
d\widetilde xdt
\le C_7\int_{  O(y)\times \{0\}}  \frac{\root\of{\tau\varphi}|f|^2 
e^{2{\tau}\alpha (t,\widetilde x',0)}}{\root\of{\lambda} } d\widetilde x'dt.
\end{equation}
By (\ref{okorok1}) and (\ref{P1Z}), we have
\begin{eqnarray}\label{P1Zz}
  \int_{\widetilde Q} ( {\tau}\lambda^2\varphi |\nabla\partial_{\widetilde x_n}
{\widetilde u}|^2 
+ {\tau}^3\lambda^4 \varphi^3 \vert \partial_{\widetilde x_n}{\widetilde u}
\vert^2)e^{2{\tau}\alpha} d\widetilde x\,dt
+ \sum_{k=1}^2\Vert \widetilde P_k(t,\widetilde x,D,\tau) 
(\partial_{\widetilde x_n}{\widetilde u}e^{\tau\alpha})\Vert^2
_{L^2({\widetilde Q})}                 \nonumber\\
\le C_8\biggl( \root\of{\tau}\int_{O(y)\times \{0\}}\root\of{\varphi} 
\frac{|f|^2 e^{2{\tau}\alpha (t,\widetilde x',0)}}{\root\of{\lambda} } 
d\widetilde x'dt   \nonumber \\
+ \int_{\widetilde\Sigma\setminus\widetilde \Sigma_0} 
\sum_{\ell=0}^1( \tau \lambda \varphi \vert\nabla\partial^\ell
_{\widetilde x_n}{\widetilde u}\vert^2
+ \tau^3 \lambda^3 \varphi^3 \vert\partial^\ell_{\widetilde x_n}
{\widetilde u}\vert^2) e^{2{\tau}\alpha} d\widetilde\Sigma
\biggr)
\end{eqnarray}
for all $ \tau\ge \tau_1(\lambda)$ and $ \lambda\ge \lambda_1$.

We set
$$
w(t,\widetilde x)={\widetilde u} e^{\tau \alpha(t,\widetilde x)}.
$$
Denote $\widetilde \Omega=O(y')\times (0,\epsilon)$.
Let $\rho = \rho(\widetilde x)\in C^2(\overline{\widetilde \Omega})$ 
satisfy 
$$
\www{n}(\widetilde x)
= \frac{\nabla\rho(\widetilde x)}{\vert\nabla \rho(\widetilde x)\vert}
\quad \mbox{on}\quad \partial\widetilde\Omega, 
$$ 
where $\www{n}(\www{x})$ is the outward unit normal vector 
to $\partial\widetilde\Omega$.
Taking the scalar product of $\widetilde P_2(t,\widetilde x,D,\tau)w$ and 
$\lambda\,\root\of{\tau}\,\root\of{\varphi}(\nabla \rho,\nabla w)$ in 
$L^2(\widetilde\Omega)$, we have
\begin{equation}
\tau^\frac 14 \root\of{\lambda}\Vert \root\of{\varphi}f e^{\tau\alpha}\Vert
_{L^2( O(y')\times\{0\})}
\le C_9(\Vert \widetilde P_2(t,\widetilde x,D,\tau)w\Vert
_{L^2(\widetilde\Omega)}
+ \lambda\,\root\of{\tau }\Vert \root\of{\varphi}\nabla w\Vert
_{L^2(\widetilde\Omega)})
\end{equation}
for all $t\in (0,T)$.

This inequality and (\ref{P1}) imply
\begin{eqnarray}\label{ZP1}
\Vert f e^{\tau\alpha}\Vert^2_{L^2(0,T;L^2(O(y')\times\{0\}))}
+ \int_{\widetilde Q} ( {\tau}\lambda\varphi |\nabla \partial_{\widetilde x_n}
{\widetilde u}|^2 
+ {\tau}^3\lambda^4 \varphi^3 \vert \partial_{\widetilde x_n}{\widetilde u}
\vert^2)
e^{2{\tau}\alpha} d\widetilde x\,dt                  \nonumber\\ 
+ \sum_{k=1}^2\Vert \widetilde P_k(t,\widetilde x,D,\tau) 
((\partial_{\widetilde x_n}{\widetilde u})e^{\tau\alpha})\Vert^2
_{L^2(\tilde Q)}
                                                      \nonumber \\
\le C_{10}\int_{\widetilde\Sigma\setminus\widetilde \Sigma_0} \sum_{\ell=0}^1
(\tau \lambda \varphi \vert\nabla\partial^\ell_{\widetilde x_n}{\widetilde u}
\vert^2
+ \tau^3 \lambda^3 \varphi^3 \vert\partial^\ell_{\widetilde x_n}{\widetilde u}
\vert^2) e^{2{\tau}\alpha} d\widetilde \Sigma.
\end{eqnarray}
Next we estimate the integral on the right-hand side of (\ref{ZP1}). 
We have
$$
\int_{\widetilde\Sigma\setminus\widetilde \Sigma_0} 
\sum_{\ell=0}^1( \tau \lambda \varphi \vert\nabla\partial^\ell_{\widetilde x_n}
{\widetilde u}\vert^2+\tau^3 \lambda^3
 \varphi^3 \vert\partial^\ell_{\widetilde x_n}{\widetilde u}\vert^2)
e^{2{\tau}\alpha}  d\widetilde \Sigma
$$
$$
= \int_{(0,T)\times O(y')\times\{\epsilon\} } 
\sum_{\ell=0}^1( \tau \lambda \varphi \vert\nabla\partial^\ell_{\widetilde x_n}
{\widetilde u}\vert^2
+ \tau^3 \lambda^3 \varphi^3 \vert\partial^\ell_{\widetilde x_n}{\widetilde u}
\vert^2) e^{2{\tau}\alpha}  d\widetilde \Sigma
$$
$$
+ \int_{(0,T)\times \partial O(y')\times (0,\epsilon) } \sum_{\ell=0}^1
( \tau \lambda \varphi \vert\nabla\partial^\ell_{\widetilde x_n}{\widetilde u}
\vert^2
+ \tau^3 \lambda^3
 \varphi^3 \vert\partial^\ell_{\widetilde x_n}{\widetilde u}\vert^2)
e^{2{\tau}\alpha} d\widetilde \Sigma=:I_1+I_2.
$$
By (\ref{klown1}), (\ref{khmer}) we can estimate $I_1$ as
\begin{eqnarray}
I_1\le \int_{(0,T)\times O(y')\times\{\epsilon\} }\sum_{\ell=0}^1
(\tau \lambda \varphi \vert\nabla\partial^\ell_{\widetilde x_n}
{\widetilde u}\vert^2
+ \tau^3 \lambda^3 \varphi^3 \vert\partial^\ell_{\widetilde x_n}{\widetilde u}
\vert^2) 
\exp\left( 2\tau\left(\alpha_I(t)-\frac{\delta}{\ell(t)} \right)\right)
d\widetilde \Sigma                                     \nonumber\\
\le \left(\int_{(0,T)\times O(y')\times\{\epsilon\}}\sum_{\ell=0}^1
(\tau \lambda \varphi \vert\nabla\partial^\ell_{\widetilde x_n}{\widetilde u}
\vert^2
+ \tau^3 \lambda^3 \varphi^3 \vert\partial^\ell_{\widetilde x_n}{\widetilde u}
\vert^2) d\widetilde \Sigma\right) e^{2{\tau}(\alpha_I(t_0)-\delta)}.
\end{eqnarray}
By (\ref{Puk}) and (\ref{klown1}), we have
\begin{eqnarray}
I_2 = \int_{(0,T)\times (\partial O(y')\setminus \mathcal G)\times(0,\epsilon)}
\sum_{\ell=0}^1( \tau \lambda \varphi \vert\nabla\partial^\ell_{\widetilde x_n}
{\widetilde u}\vert^2
+ \tau^3 \lambda^3 \varphi^3 \vert\partial^\ell_{\widetilde x_n}{\widetilde u}
\vert^2)e^{2{\tau}\alpha} d\widetilde \Sigma                 \nonumber\\
\le \int_{(0,T)\times (\partial O(y')\setminus F(\mathcal G))
\times (0,\epsilon)} 
\sum_{\ell=0}^1( \tau \lambda \varphi \vert\nabla\partial^\ell
_{\widetilde x_n}{\widetilde u}\vert^2
+ \tau^3 \lambda^3
 \varphi^3 \vert\partial^\ell_{\widetilde x_n}{\widetilde u}\vert^2)\times\nonumber\\\times
\exp\left( \frac{2\tau\left( e^{\lambda K} - e^{2\lambda ||\psi||
_{C^0(\overline{\tilde \Omega})}}\right)}{\ell(t)}\right) 
d\widetilde x'dt                                         \nonumber\\ 
\le \left(\int_{(0,T)\times (\partial O(y')
\setminus F(\mathcal G))\times(0,\epsilon)}
\sum_{\ell=0}^1( \tau \lambda \varphi \vert\nabla\partial^\ell_{\widetilde x_n}
{\widetilde u}\vert^2
+ \tau^3 \lambda^3
 \varphi^3 \vert\partial^\ell_{\widetilde x_n}{\widetilde u}\vert^2)
d\widetilde x'dt \right)        \nonumber \\
\times e^{2{\tau} (e^{\lambda K}
- e^{2\lambda ||\psi||_{C^0(\overline{\widetilde \Omega})}}) }.
\end{eqnarray}

In view of (\ref{ZP1}), there exists a constant $\delta>0$ such that
\begin{eqnarray}\label{zZP1}
\Vert f e^{\tau\alpha_I}\Vert^2_{L^2(0,T;L^2(O(y')\times\{0\}))}
+ \int_{\widetilde Q} ( {\tau}\lambda^2\varphi |\nabla 
\partial_{\widetilde x_n}{\widetilde u}|^2 
+ {\tau}^3 \lambda^4\varphi^3 \vert \partial_{\widetilde x_n}{\widetilde u}
\vert^2) e^{2{\tau}\alpha} d\widetilde x\,dt                      \nonumber\\
+ \sum_{k=1}^2\Vert \widetilde P_k(t,\widetilde x,D,\tau)
 ((\partial_{\widetilde x_n}{\widetilde u})e^{\tau\alpha})\Vert^2
_{L^2(\widetilde Q)}                          \nonumber\\ 
\le C_{11}\biggl(
\int_{(0,T)\times O(y')\times \{\epsilon\}}
\sum_{\ell=0}^1( \tau \lambda \varphi \vert\nabla\partial^\ell_{\widetilde x_n}
{\widetilde u}\vert^2
+ \tau^3 \lambda^3 \varphi^3 \vert\partial^\ell_{\widetilde x_n}{\widetilde u}
\vert^2)
\exp\left( 2\tau\left( \alpha_I(t)-\frac{\delta}{\ell(t)} \right)\right)
  d\widetilde \Sigma                                          \nonumber\\
+ \left( \int_{(0,T)\times \partial O(y')\times (0,\epsilon) } \sum_{\ell=0}^1
( \tau \lambda \varphi \vert\nabla\partial^\ell_{\widetilde x_n}{\widetilde u}
\vert^2 + \tau^3 \lambda^3 \varphi^3 \vert\partial^\ell_{\widetilde x_n}
{\widetilde u}\vert^2) d\widetilde x'dt \right)
e^{2{\tau}(e^{\lambda K} - e^{2\lambda ||\psi||
_{C^0(\overline{\widetilde \Omega})}}) }\biggr).
\end{eqnarray}
By (\ref{klown}) and (\ref{klown1}), there exist constants $C_{12}>0$
and $m>1$ such that
\begin{equation}\label{kokom}
\exp\left( 2\tau\left(\alpha_I(t_0)-\frac{\delta}{m}\right)\right)
\le C_{12}\Vert f e^{\tau\alpha_I}\Vert^2_{L^2(0,T;L^2(O(y')\times\{0\}))}.
\end{equation}

By (\ref{kokom}) and (\ref{zZP1}), there exists a constant $C_{13}>0$ 
such that
\begin{eqnarray}\label{KzzZP1}
\Vert f e^{\tau\alpha_I}\Vert^2_{L^2(0,T;L^2(O(y')\times\{0\}))}
+ \int_{\widetilde Q} ( {\tau}\lambda^2\varphi |\nabla 
\partial_{\widetilde x_n}{\widetilde u}|^2 
+ {\tau}^3 \lambda^4\varphi^3 \vert \partial_{\widetilde x_n}{\widetilde u}
\vert^2)e^{2{\tau}\alpha}   d\widetilde x\,dt                 \nonumber\\  
\le C_{13} e^{2{\tau} (e^{\lambda K} - e^{2\lambda ||\psi||
_{C^0(\overline{\widetilde \Omega})}}) }\int_{(0,T)\times \partial O(y')\times (0,\epsilon) } 
\sum_{\ell=0}^1( \tau \lambda \varphi \vert\nabla\partial^\ell
_{\widetilde x_n}{\widetilde u}\vert^2
+ \tau^3 \lambda^3 \varphi^3 \vert\partial^\ell_{\widetilde x_n}{\widetilde u}
\vert^2) d\widetilde x'dt .
\end{eqnarray}

Since $\widetilde \psi(y')>0$, there exist an open ball $B(\delta_2,y')\subset 
\Bbb R^{n-1}$ and a constant $\delta_3>0$ such that
$$
\alpha(t,\widetilde x',0)
\ge (e^{\lambda K} - e^{2\lambda ||\psi||_{C^0(\overline{\widetilde \Omega})}})
-\delta_3 \quad \mbox{for all $t\in [t_0-\delta_1,t_0+\delta_1]$ and 
$\widetilde x'\in B(\delta_2,y')$}.
$$
Hence, this inequality and (\ref{KzzZP1}) yield
\begin{eqnarray}\label{pauk}
\Vert f e^{2\tau (e^{\lambda K}
- e^{2\lambda ||\psi||_{C^0(\overline{\widetilde \Omega})}})-\delta_3)}
\Vert^2_{L^2(t_0-\delta_1,t_0+\delta_1;L^2(O(y')\times\{0\}))}    \\
\le\Vert f e^{\tau\alpha_I}\Vert^2_{L^2(t_0-\delta_1,t_0+\delta_1;L^2(O(y')
\times\{0\}))}
\le \Vert f e^{\tau\alpha_I}\Vert^2_{L^2(0,T;L^2(O(y')\times\{0\}))}
                                           \nonumber\\
\le C_{14} e^{2{\tau} (e^{\lambda K} - e^{2\lambda ||\psi||
_{C^0(\overline{\widetilde \Omega})}}) }
\int_{(0,T)\times \partial O(y')\times (0,\epsilon) } \sum_{\ell=0}^1
(\tau \lambda \varphi \vert\nabla\partial^\ell_{\widetilde x_n}{\widetilde u}
\vert^2
+ \tau^3 \lambda^3
 \varphi^3 \vert\partial^\ell_{\widetilde x_n}{\widetilde u}\vert^2)  
d\widetilde x'dt .\nonumber
\end{eqnarray}
Passing in (\ref{pauk}) as $\tau$ goes to $+\infty$, we can find 
a neighborhood $U$ of $y'$ such that $f = 0$ in $(0,T)\times U$.
$\blacksquare$

{\bf Proof of Theorem \ref{lox}.} Our proof is by contradiction. Suppose that
there exist $(u_1,f_1)$ and $(u_2,f_2)$ which satisfy  the conditions of 
Theorem \ref{lox} and (\ref{gavno}). Setting $u=u_1-u_2$ and $f=f_1-f_2$,
we have
\begin{equation}\label{KSsok1}
P(t,x,D)u=R(t,x) f(t,x')\quad \mbox{in}\quad Q,
\end{equation}
\begin{equation}\label{KSSsok2}
u\vert_{\Sigma}=0, \quad \partial_\nu u_1=\partial_\nu u_2\quad \mbox{on}\,\, 
(0,T)\times \Gamma, \quad f_1=f_2 \quad \mbox{on}\,\, (0,T)\times \Gamma_1.
\end{equation}
Let $z_1 \in S_1(y',y_n)\cap S_2(y',y_n)$.
By Lemma \ref{lox1}, we see that $f(t,x') = 0$ in $(0,T)\times U(z_1)$,
where $U(z_1)$ is a neighborhood of $z_1$.

Let $k_0\in \N$ be the maximal number satisfying: there exist a point 
$z_{k_0-1}\in S_{k_0}(y',y_n)$ and a neighborhood  $U(z_{k_0-1})$ of 
$z_{k_0-1}$ such that $f(t,x') = 0$ in $(0,T)\times U(z_{k_0-1})$.
If $k_0=K(y')$, then Lemma \ref{lox1} implies that there exists a neighborhood of $W$ 
of $(y',y_n)$ such that $f=0$ in $(0,T)\times W$.
Otherwise, there exist a point $z_{k_0}\in S_{k_0}(y',y_n)\cap 
S_{k_0+1}(y',y_n)$ and a neighborhood $U(z_{k_0})$ of $z_{k_0}$ such that 
$f(t,x')=0$ in $(0,T)\times U(z_{k_0})$, and so we reached a contradiction 
to the fact that $k_0$ is maximal.
$\blacksquare$

{\bf Proof of Theorem \ref{Plox}.} Let $u_1$ and $u_2$ satisfy (\ref{PKSsok1}) 
and (\ref{MKSsok1}) respectively. We set $u=u_1-u_2$ and $f=c_1-c_2.$
By (\ref{PKSsok1}), (\ref{MKSsok1}) and (\ref{KKSSsok2}), we have
\begin{equation}\label{mk}
\partial_tu-\sum_{i,j=1}^n \partial_{x_i}(a_{ij}(x)\partial_{x_j} u)
+\sum_{i=1}^n b_{i}(x)\partial_{x_i} u+c_1(t,x') u=-u_2c\quad \mbox{in}\quad Q
\end{equation}
and
\begin{equation}\label{mk1}
u\vert_\Sigma=0.
\end{equation}
By (\ref{gavno}), we have
$$
\partial_\nu u\vert_{(0,T)\times \Gamma}=0,\quad 
f\vert_{(0,T) \times \Gamma_1}=0.
$$
We set $R(t,x)=-u_2(t,x).$ By (\ref{gavno1})
there exists a constant $\beta>0$ such that
$$
\vert R(t,x)\vert>\beta>0\quad \mbox{on}\quad  \ooo{Q}.
$$ 
Applying Theorem \ref{lox} to (\ref{mk}) and (\ref{mk1}),
we complete the proof of Theorem \ref{Plox}.
$\blacksquare$

\section{Proof of Theorem \ref{XPlox}.}
{\bf Proof of Theorem \ref{XPlox}.}
Let $u_1$ and $u_2$ satisfy (\ref{XPKSsok1}) and (\ref{XMKSsok1}) 
respectively with $c_1$ and $c_2$. We set $u=u_1-u_2,$  $f=c_1-c_2$ and
$$
A(x,D)u:=-\sum_{i,j=1}^n \partial_{x_i}(a_{ij}(x)\partial_{x_j} u)
+\sum_{i=1}^n b_{i}(x)\partial_{x_i} u+c_1(x') u.
$$
By (\ref{XMKSsok1}) we have
\begin{equation}
\sqrt{-1}\partial_tu+A(x,D)u=-u_2f\quad \mbox{in}\quad Q,
\end{equation}
\begin{equation}
u\vert_\Sigma=0.
\end{equation}
By (\ref{Pgavno})
$$
\partial_\nu u\vert_{(0,T)\times \Gamma}=0\quad\mbox{and} \quad 
f\vert_{(0,T)\times \Gamma_1}=0.
$$
We set $R(t,x)=-u_2(t,x)$.
By (\ref{Pgavno1}), there exists a constant $\beta>0$ such that
$$
\vert R(t,x)\vert>\beta>0\quad \mbox{on}\quad  \ooo Q.
$$

Let $\mbox{\bf G}=\mbox{Int}\,\{x'\in \mathcal O\vert f(x')=0\}$.
Suppose that $\mathcal O\setminus \mbox{\bf G}\ne\{\emptyset\}$.
Then there exists a point $y'\in \partial\mbox{\bf G}\cap \mathcal O$.
By Condition 2 one can find  $y_n\in \Bbb R^1$ such that $ (y',y_n)\in \Gamma.$ Hence there exists a neighborhood $U$ of $y$ such that $U\cap \Gamma$ is 
a graph of some smooth  function $g(x').$

Let $\mathcal K(x_0,\tau)$ be a solution to the boundary controllability problem for the one dimensional  Schr{\"o}dinger equation:
\begin{eqnarray}\label{HUKK}
\sqrt{-1}\partial_{\tau}\mathcal K-\partial^2_{x_0}\mathcal K =0\quad  
\mbox{in} \quad  x_0\in (0,1), \,\tau\in(0,T),\nonumber\\
\quad \mathcal K(0,\tau)=\psi(\tau)\,\,\mbox{on}\,\, (0,T),
\quad \mathcal K(\cdot,0)=\mathcal K(\cdot,T)=0.
\end{eqnarray}

A control is located on the part of the boundary $t=1:$
$$
\mathcal K(1,\tau)= v(\tau), \quad 0<\tau<T.
$$ 
For given $\psi\in C^\infty_0(0,T)$, the existence of such a smooth solution 
$K$ to problem (\ref{HUKK}) is proved in \cite{MRR}.

Consider 
$$
w(x_0,x)=\int_{0}^T\mathcal K(x_0,\tau)u(\tau,x)d\tau, \quad 
M(x_0,x)=\int_{0}^T\mathcal K(x_0,\tau)R(\tau,x)d\tau.
$$
Observe that $w \in H^2((0,1)\times \Omega)$ satisfies 
\begin{equation}\label{Ssimas}
L(x,D) w:=-\partial_{x_0}^2w+A(x,D)w=M(x_0,x)f(x') \quad \mbox{in}\quad
  (0,1)\times  \Omega,
\end{equation}
\begin{equation}
\quad w=0\quad \mbox{on $(0,1)\times \partial\Omega$},\quad  
\partial_\nu w=0\quad \mbox{on $(0,1)\times \Gamma$}, 
\quad f=0\quad \mbox{on $\Gamma_1$}.
\end{equation}
We choose a function $\psi \in C^{\infty}_0(0,T)$ such that
$$
 M(0,y) =\int^T_0 \psi(\tau)R(\tau,y) d\tau \ne 0.
$$ 
Then there exist constants $b\in (0,1)$ and $\beta_1>0$, and 
a neighborhood $U$ of $y$ such that
\begin{equation}
\vert M(x_0,x)\vert>\beta_1>0 \quad \mbox{ for all }\,\,(x_0,x) \in \quad [0,b]\times 
(U\cap \Omega).
\end{equation}

The rest of the proof is similar to the proof of Lemma \ref{lox1} as follows.
We make a change of the variables $x \mapsto \www{x}$ in a neighborhood  
$\Psi$ of the surface $S(y',y_n)\cap U:$
$$
\widetilde x=(\widetilde x',\widetilde x_n)=F(x)=(x', x_n-g(x')).
$$ 
An unknown function $f$ is not affected by this change of variables. 
By $\widetilde w$ and $\widetilde M$ we denote the functions $w$ and $M$ 
in the new variables.

On the other hand,  there exists a constant $\epsilon>0$  such that
$$
\widetilde \Omega=O(y')\times(0,\epsilon)\subset F(\Psi\cap \Omega).
$$
The change of coordinates transforms the elliptic operator $L$ into an
elliptic operator $\widetilde L$:
\begin{equation}\label{Kfocus}
\widetilde L(\widetilde x,D)\widetilde w
= -\partial^2_{x_0}\widetilde w-\sum_{i,j=1}^n \partial_{\widetilde x_i}
(\widetilde a_{ij}(\widetilde x)\partial_{\widetilde x_j}\widetilde  w)
\end{equation}
$$
+ \sum_{i=1}^n \tilde b_{i}(\widetilde x)\partial_{\widetilde x_i} 
\widetilde w
+ \widetilde c(x') \widetilde w
= \widetilde M(x_0,\widetilde x) f(\widetilde x')\quad \mbox{in}\quad (0,b)
\times F(\Psi\cap \Omega),
$$
\begin{equation}\label{KSsok2}
\partial_{\widetilde x_n}\widetilde w\vert_{(0,b)\times O(y')\times\{0\}}
= \widetilde w\vert_{(0,b)\times O(y')\times\{0\}}=0.
\end{equation}

We set $\widetilde a_{00}=1$ and $\widetilde a_{0,j}=0$ for all 
$j\in \{1,\dots,n\}$.
By (\ref{sok4})-(\ref{sok6}), the coefficients $\widetilde a_{ij}$,
$\widetilde b_j$ and $\widetilde c$ satisfy

\begin{equation}\label{Ksok4}
\widetilde a_{ij}=\widetilde a_{j i}\quad \mbox{for all
$i,j\in \{0,\dots,n\}$ and all $\widetilde x\in F(\Psi\cap \Omega)$}
\end{equation}
and
\begin{eqnarray} 
\widetilde a_{ij},\, 
\partial_{\widetilde x_n}\widetilde a_{ij} \in 
C^{1}(\overline{ F(\Psi\cap \Omega)}), \quad 
\widetilde b_{i},\partial_{\widetilde x_n}\widetilde b_i\in 
C^{0}(\overline{ F(\Psi\cap \Omega)}),                   \nonumber\\
\quad \widetilde c,\partial_{\widetilde x_n} \widetilde c\in  C^{0}
(\overline{ F(\Psi\cap \Omega)})\quad \mbox{for all $i,j\in \{0,\dots,n\}$}
\end{eqnarray}
and we can find constants $\beta>0$ and $\www{\beta}_1 > 0$ such that
\begin{equation}\label{Ksok5}
\sum_{i,j=0}^n \widetilde a_{ij}(\widetilde x)\eta_{i}\eta_j
\ge \beta\vert \eta\vert^2\quad \mbox{for all $(\widetilde x,\eta)
\in F(\Psi\cap \Omega)\times \Bbb R^n$}
\end{equation}
and 
\begin{equation}\label{Ksok6}
\widetilde M, \, \partial_{\widetilde x_n}\widetilde M\in 
C^{2}(\overline{(0,b)\times F(\Psi\cap \Omega)}), \quad 
\vert \widetilde M(x_0,\tilde x)\vert >\widetilde \beta_1>0\quad 
\mbox{on $\overline{(0,b)\times F(\Psi\cap \Omega)}$}.
\end{equation}
Let us introduce a weight function
\begin{equation}\label{Ksous}
\psi(x_0,\widetilde x)= \widetilde \psi(\widetilde x')-N\widetilde x_n^2
-\widetilde N(\widetilde x_0-a)^2 +K.
\end{equation}
Here we can choose positive constants $N$, $\widetilde N$, $K$ and 
$a \in (0,b)$ such that 
\begin{equation}\label{Kkkm}
\inf_{\widetilde x'\in \mathcal O(y')}\psi(x_0,\widetilde x',0)
> \sup_{\widetilde x'\in \mathcal O(y')}\psi(x_0,\widetilde x',\epsilon)\quad 
\mbox{for all $x_0\in [0,b]$},
\end{equation} 
\begin{equation}\label{KPuk1} 
\psi(x_0,\widetilde x)>0\quad \mbox{on}\,\,\overline{(0,b)\times\widetilde 
\Omega}.
\end{equation}
Let $\widetilde \psi\in C^2(\overline{O(y')})$ satisfy 
\begin{eqnarray} \label{KPuk}
\nabla\widetilde\psi\ne 0\quad\mbox{on}\,\, \overline{O(y')},\quad  
\partial_n \widetilde \psi<0\quad \mbox{on} \quad \partial O(y'),\nonumber\\ \quad 
\widetilde \psi\vert_{\partial O(y')\setminus \widetilde{\mathcal G}}=0, 
\quad \widetilde\psi >0 \quad\mbox{on $O(y')\setminus 
\widetilde{\mathcal G}$}, 
\end{eqnarray} 
where $\widetilde{\mathcal G}=F(\mbox{\bf G}).$
The existence of such a function is established in \cite{Ima1995}.

In terms of $\psi$, we introduce a function  $\varphi(x_0,\widetilde x)$:
\begin{equation}\label{KDD1.2} 
\varphi(x_0, \widetilde x) = e^{\lambda \psi(x_0,\widetilde x)},
\end{equation}
where $\lambda$ is a large positive parameter.

By (\ref{Kkkm}), we have
\begin{equation}\label{Kkhmer}
\varphi_I(x_0)=\inf_{\widetilde x'\in \mathcal O(y')}
\varphi(x_0,\widetilde x',0)> \sup_{\widetilde x'\in \mathcal O(y')}
\varphi(x_0,\widetilde x',\epsilon)=\varphi_M(x_0) \quad 
\mbox{for all $x_0\in [0,b]$}
\end{equation}
and there exist constants $\delta_1>0, \delta_2>0$ and $\delta_3>0$  such that
\begin{equation}\label{soska} 
\varphi(x_0,\widetilde x',0) >e^{\lambda K}+\delta_2 \quad 
\mbox{for all $(x_0,\widetilde x')\in [a-\delta_1,a+\delta_1]\times 
B(y',\delta_3)$}.
\end{equation}

We assume that $\mu:= \mu(x_0)\in C_0^\infty(0,b)$ satisfies 
\begin{equation}\label{proton}
\mu\vert_{[\delta_4,b-\delta_4]}=1, \quad \varphi(a,\widetilde x)
> \mbox{max}\{\Vert \varphi\Vert_{C^0([b-\delta_4,b]\times \ooo\Omega)}, \,
\Vert \varphi\Vert_{C^0([0,\delta_4]\times \ooo\Omega)}\},
\end{equation}
where $\delta_4>0$ is some constant.
For each sufficiently large parameter $\widetilde N$, such a function $\mu$ 
exists.

The function $\mu \widetilde w$ satisfies
\begin{equation}\label{KKfocus}
\widetilde L(\widetilde x,D)(\mu \widetilde w)
= \mu\widetilde M(x_0,\widetilde x) f(\widetilde x')
+ [\mu,\widetilde L]\widetilde w\quad \mbox{in}\quad (0,b)\times 
F(\Psi\cap \Omega),
\end{equation}
\begin{equation}\label{KKSsok2}
\partial_{\widetilde x_n}(\mu\widetilde w)\vert_{(0,b)\times O(y')\times\{0\}}
= \mu\widetilde w\vert_{(0,b)\times O(y')\times\{0\}}=0.
\end{equation}
We consider a boundary value problem
\begin{equation}\label{KLD1.4}
 \widetilde L(\widetilde x,D)z =  g \quad \mbox{in}\quad G\triangleq(0,b)
\times O(y')\times (0,\epsilon),
\end{equation}
\begin{equation}\label{KLD1.5}
z \big|_{\widetilde S_0} = 0,\quad z(0, \cdot) = z(b,\cdot)=0,
\end{equation}
where
$\widetilde S_0=(0,T)\times O(y')\times \{0\}$.

Henceforth we set  
$$
\gamma=\sum_{j=1}^n\Vert \widetilde b_j\Vert_{L^\infty(G)}
+ \Vert \widetilde c\Vert_{L^\infty(G)}, \quad
\widetilde \nabla=(\partial_{\widetilde x_0},\partial_{\widetilde x_1},\dots,
\partial_{\widetilde x_n}).
$$

In order to formulate the next lemma, we introduce two operators
\begin{equation}\label{sosna}
\left\{ \begin{array}{rl}
& \widetilde L_1(x_0,\widetilde x,D,\tau)z
:= 2\tau\varphi\lambda \widetilde a(\widetilde x,\widetilde\nabla \psi,
\widetilde \nabla z)\\
+ &\tau\lambda^2\varphi \widetilde a(\widetilde x,\widetilde\nabla\psi,
\widetilde\nabla\psi)z
+ \tau\varphi\lambda \sum_{i,j=0}^n(\partial_{\widetilde x_i}\widetilde a_{ij})
(\partial_{\widetilde x_j}\psi) z, \\
& \widetilde L_2(x_0,\widetilde x,D,\tau)z
:= -\partial^2_{x_0} z- \sum_{i,j=1}^n \widetilde a_{ij} 
\partial^2_{\widetilde x_i\widetilde x_j} z -
\lambda^2 {\tau}^2 \varphi^2 \widetilde a(\widetilde x, \widetilde\nabla\psi, 
\widetilde\nabla\psi) z,
\end{array}\right.
\end{equation}
where $\widetilde a(\widetilde x,\eta,\eta)
= \sum_{i,j=0}^n \widetilde a_{ij}(\widetilde x)\eta_i\eta_j$.

Then
\begin{lemma}\label{KD1.2} 
{\it Let (\ref{Ksok4}) and (\ref{Ksok5}) hold true,
and the function $\varphi$ be defined as in (\ref{Ksous}) and (\ref{KDD1.2}), 
where $\psi$ satisfies (\ref{KPuk1}) and (\ref{KPuk}). 
Let $z\in H^{2}(G)$ and $g\in L^2(G)$.
Then there exists a constant $\widehat\lambda(\gamma) = \widehat\lambda> 0$ 
such that 
for an arbitrary $\lambda \ge \widehat{\lambda}$ there exists
a constant ${\tau}_0(\lambda,\gamma)>0$ such that any solutions 
of problem (\ref{KLD1.4}) - (\ref{KLD1.5})
satisfy the following inequality:
\begin{eqnarray}\label{KDD1.6}
\int_{G} ( {\tau}\lambda^2\varphi |\widetilde\nabla z|^2 
+ {\tau}^3\lambda^4 \varphi^3 \vert z\vert^2)
e^{2{\tau}\varphi}   d\widetilde x\,dx_0
+ \sum_{k=1}^2\Vert \widetilde L_k(x_0,\widetilde x,D,\tau) 
(ze^{\tau\varphi})\Vert^2_{L^2(G)}\nonumber\\
 \le C_1\left( \int_{G} |g|^2 e^{2{\tau}\varphi}  d\widetilde xdx_0
+ \int_{\partial G\setminus\widetilde {S}_0} 
(\tau \lambda\varphi \vert\widetilde \nabla z\vert^2+\tau^3 \lambda^3\varphi^3 
\vert z\vert^2)e^{2{\tau}\varphi}  d\widetilde {S}\right),
\end{eqnarray}
for each ${\tau} \ge \tau_0(\lambda,\gamma)$.
Here the constant $C_1$ is independent of $\tau$ and $\lambda$.
}
\end{lemma}

Applying to the  system (\ref{Kfocus}) Carleman estimate (\ref{KDD1.6}), 
we obtain
\begin{eqnarray}\label{KP1}
\int_{G} ( {\tau}\lambda^2\varphi |\widetilde \nabla\widetilde  w|^2 
+ {\tau}^3 \lambda^4\varphi^3 \vert \widetilde w\vert^2)
e^{2{\tau}\varphi}   d\widetilde x\,dx_0
+ \sum_{k=1}^2\Vert\widetilde L_k(x_0,\widetilde x,D,\tau) 
(\widetilde w e^{\tau\varphi})\Vert^2_{L^2({G})}             \nonumber\\
\le C_2\left( \int_{G} (|f|^2 +\vert [\mu,\widetilde L]\widetilde w\vert^2)
e^{2{\tau}\varphi}  d\widetilde x dx_0 
+ \int_{\partial G\setminus\widetilde {S}_0} 
(\tau \lambda\varphi \vert\widetilde \nabla \widetilde w\vert^2
+ \tau^3 \lambda^3\varphi^3 \vert \widetilde w\vert^2)
e^{2{\tau}\varphi}  d\widetilde{S}\right)
\end{eqnarray}
for all $ \tau\ge \tau_0(\lambda)$ and $ \lambda\ge \lambda_0$.

Next we differentiate both sides of equation (\ref{Kfocus}) 
with respect to $\widetilde x_n$:
\begin{equation}\label{Kfocus1}
\widetilde L(\widetilde x,D)\partial_{\widetilde x_n}\widetilde w
= [\partial_{\widetilde x_n},\widetilde L] \widetilde w
+ \partial_{\widetilde x_n}\widetilde M(x_0,\tilde x) f(\widetilde x')
\quad \mbox{in}\quad G,
\end{equation}
\begin{equation}\label{KfSsok2}
\partial_{\widetilde x_n}\widetilde w\vert_{(0,b)\times O(y')\times\{0\}}=0.
\end{equation}
Applying to the system (\ref{Kfocus1}) Carleman estimate (\ref{KDD1.6}),
we obtain
\begin{eqnarray}\label{KP11}
\int_{G} ( {\tau}\lambda^2\varphi |\widetilde \nabla \partial_{\widetilde x_n}
\widetilde w|^2 
+ {\tau}^3 \lambda^4\varphi^3 \vert \partial_{\widetilde x_n}\widetilde w
\vert^2)e^{2{\tau}\varphi}  d\widetilde x\,dx_0\nonumber\\
+ \sum_{k=1}^2\Vert \widetilde L_k(x_0,\widetilde x,D,\tau) 
((\partial_{\widetilde x_n}\widetilde w)e^{\tau\varphi})\Vert^2_{L^2(G)}
                                   \nonumber\\
\le C_3\biggl( \int_{G} \{ 
(|f|^2 +\vert \partial_{\widetilde x_n}[\mu,\widetilde L]\widetilde w\vert^2) 
e^{2{\tau}\varphi}  
+ \sum_{\vert\gamma\vert\le 2}\vert \partial^\gamma_{\widetilde x} 
\widetilde w\vert^2 e^{2\tau\varphi} \} d\widetilde xdx_0   \nonumber\\
+ \int_{\partial G\setminus\widetilde {S}_0} 
(\tau \lambda \varphi \vert\widetilde \nabla\partial_{\widetilde x_n}
\widetilde w\vert^2 + \tau^3 \lambda^3 \varphi^3 \vert\partial
_{\widetilde x_n}\widetilde w\vert^2)
e^{2{\tau}\varphi}  d\widetilde {S}\biggr)
\end{eqnarray}
for all $ \tau\ge \tau_1(\lambda)$ and $ \lambda\ge \lambda_1$.

From (\ref{KP1}) and (\ref{KP11}), we have
\begin{eqnarray}\label{KP1Z}
\int_{G} ( {\tau}\lambda^2\varphi |\widetilde \nabla \partial_{\widetilde x_n}
\widetilde w|^2 + {\tau}^3\lambda^4 \varphi^3 \vert 
\partial_{\widetilde x_n}\widetilde w\vert^2)e^{2{\tau}\varphi} 
d\widetilde x\,dx_0\nonumber\\
+ \sum_{k=1}^2\Vert \widetilde L_k(x_0,\widetilde x,D,\tau) 
((\partial_{\widetilde x_n}\widetilde w)e^{\tau\varphi})\Vert^2_{L^2({G})}
                                    \nonumber\\
\le C_4\biggl( \int_{G}(\varphi \tau (|f|^2
+ \vert [\mu,\widetilde L]\widetilde w\vert^2)
+ \vert\partial_{\widetilde x_n}[\mu,\widetilde L]\widetilde w\vert^2 )
e^{2{\tau}\varphi}  d\widetilde xdx_0              \nonumber \\
+ \int_{\widetilde{S}\setminus\widetilde {S}_0}\sum_{\ell=0}^1 
(\tau \lambda \varphi \vert\partial^\ell_{\widetilde x_n}\widetilde w\vert^2
+ \tau^3 \lambda^3 \varphi^3 \vert\partial^\ell_{\widetilde x_n}\widetilde w
\vert^2) e^{2{\tau}\varphi}  d\widetilde{S}\biggr)
\end{eqnarray}
for all $ \tau\ge \tau_2(\lambda)$ and $ \lambda\ge \lambda_2.$
By (\ref{Ksous}), we have
$$
\int_{G}\varphi |f|^2 e^{2{\tau}\varphi (x_0,\widetilde x)}  d\widetilde xdx_0
$$
$$
= \int_{(0,b)\times O(y')} e^{\lambda\widetilde \psi(\widetilde x')
-\lambda\widetilde N(x_0-a)^2 + \lambda K } 
|f|^2  \int_0^\epsilon  e^{-\lambda N\widetilde x_n^2}
e^{2\tau  e^{\lambda\widetilde \psi(\widetilde x')-\lambda\widetilde N(x_0-a)^2
+\lambda K} e^{-\lambda N \widetilde x_n^2}}d\widetilde x_n d\widetilde x'dx_0=
$$
$$
\int_{(0,b)\times O(y') } e^{\lambda\widetilde \psi(\tilde x')
-\lambda\widetilde N(x_0-a)^2+\lambda K } |f|^2 e^{2{\tau}\varphi (x_0,\widetilde x',0)}
$$
$$
\times  \left(\int_0^{\epsilon\sqrt{\lambda N}}  
e^{-\widetilde x_n^2}e^{2\tau e^{\lambda\widetilde \psi(\tilde x')
-\lambda\widetilde N(x_0-a)^2+\lambda K } 
(e^{-\tilde x_n^2}-1)}d\widetilde x_n
\right) \frac{1}{\sqrt{\lambda}} d\widetilde x'dx_0.
$$
Using the Morse lemma (see e.g. \cite{Her}) one can find a coordinates $\tilde x_n=r(y_n)$ such that $(e^{-\tilde x_n^2}-1)\circ r(y_n)=-y_n^2.$
Making the corresponding change of variables in the integral  we have 
\begin{eqnarray}\label{Kokorok}
\int_0^{\epsilon\sqrt{\lambda N}}  
e^{-\widetilde x_n^2}e^{2\tau e^{\lambda\widetilde \psi(\tilde x')
-\lambda\widetilde N(x_0-a)^2+\lambda K } 
(e^{-\tilde x_n^2}-1)}d\widetilde x_n\le C_5\int_{-\infty}^{+\infty}  
e^{-2\tau e^{\lambda\widetilde \psi(\tilde x')
-\lambda\widetilde N(x_0-a)^2+\lambda K } 
y_n^2}d y_n \nonumber\\\le \frac{C_6}{\root\of{\tau\varphi(x_0,\tilde x',0)}}.
\end{eqnarray}

Using (\ref{Kokorok}), we obtain
\begin{equation}\label{Kokorok1}
\int_{G}\tau \varphi |f|^2 e^{2{\tau}\varphi (x_0,\widetilde x)}  
d\widetilde xdx_0
\le C_7\int_{\widetilde  O(y)\times \{0\}} \root\of{\tau\varphi} 
\frac{|f|^2 e^{2{\tau}\varphi (x_0,\widetilde x',0)}}{\root\of{\lambda} } 
d\widetilde x'dx_0.
\end{equation}
By (\ref{Kokorok1}) and (\ref{KP1Z}), we have
\begin{eqnarray}\label{KP1Zz}
\int_{G} ( {\tau}\lambda^2\varphi |\widetilde \nabla \partial_{\widetilde x_n}
\widetilde w|^2 
+ {\tau}^3\lambda^4 \varphi^3 \vert \partial_{\widetilde x_n}\widetilde w
\vert^2)
e^{2{\tau}\varphi} d\widetilde x\,dx_0\nonumber\\
+ \sum_{k=1}^2\Vert \widetilde L_k(x_0,\widetilde x,D,\tau) 
((\partial_{\widetilde x_n}\widetilde w)
e^{\tau\varphi})\Vert^2_{L^2({G})}\nonumber\\
\le C_8\biggl( \root\of{\tau}\int_{\widetilde O(y)\times \{0\}}
\root\of{\varphi} \frac{|f|^2 e^{2{\tau}\varphi (x_0,\widetilde x',0)}}
{\root\of{\lambda} } d\widetilde x'dx_0 \nonumber\\
+ \int_{G}(\varphi \tau \vert [\mu,\widetilde L]\widetilde w\vert^2
+ \vert\partial_{\widetilde x_n}[\mu,\widetilde L]\widetilde w\vert^2 )
e^{2{\tau}\varphi}  d\widetilde xdx_0                  \nonumber\\
+ \int_{\partial G\setminus\widetilde {S}_0} 
\sum_{\ell=0}^1( \tau \lambda \varphi \vert\widetilde \nabla\partial^\ell
_{\widetilde x_n}\widetilde w\vert^2
+ \tau^3 \lambda^3 \varphi^3 \vert\partial^\ell_{\widetilde x_n}
\widetilde w\vert^2)e^{2{\tau}\varphi}  d\widetilde{S}\biggr)
\end{eqnarray}
for all $\tau\ge \tau_1(\lambda)$ and $ \lambda\ge \lambda_1$.

We set
$$
\mbox{\bf w}(x_0,\widetilde x)=w e^{\tau \varphi(x_0,\widetilde x)}.
$$
Denote $\widetilde \Omega=O(y')\times (0,\epsilon).$
Let  $\rho(\widetilde x)\in C^2(\overline{\widetilde \Omega})$ satisfy
$$
n(\widetilde x)=\frac{\nabla\rho(\widetilde x)}{\vert\nabla \rho(\widetilde x)
\vert} \quad \mbox{on $\partial\widetilde\Omega$},
$$ 
where $n(\widetilde{x})$ denotes the outward unit normal vector to 
$\partial\widetilde\Omega.$
Taking the scalar product of the functions 
$\widetilde L_2(x_0,\widetilde x,D,\tau)w$ and 
$\root\of{\tau}\lambda\,\,\root\of{\varphi}(\widetilde \nabla \rho,
\widetilde \nabla \mbox{\bf w})$ in $L^2(\widetilde\Omega)$, we have
\begin{equation}
\tau^\frac 14 \root\of{\lambda}\Vert \root\of{\varphi}f e^{\tau\varphi}\Vert
_{L^2( O(y')\times\{0\})}
\le C_9(\Vert \widetilde L_2(x_0,\widetilde x,D,\tau)\mbox{\bf w}\Vert
_{L^2(\widetilde\Omega)}
+  \root\of{\tau }\lambda\Vert \root\of{\varphi}\widetilde \nabla\mbox{\bf w}
\Vert_{L^2(\widetilde\Omega)})
\end{equation} or all $x_0\in (0,b)$.

This inequality and (\ref{KP1}) imply
\begin{eqnarray}\label{KZP1}
\Vert f e^{\tau\varphi}\Vert^2_{L^2(0,b;L^2(O(y')\times\{0\}))}
+ \int_{G} ( {\tau}\lambda\varphi |\widetilde \nabla \partial_{\widetilde x_n}
\widetilde w|^2 + {\tau}^3\lambda^4 \varphi^3 \vert \partial_{\widetilde x_n}
\widetilde w\vert^2)
e^{2{\tau}\varphi} d\widetilde x\,dx_0                       \nonumber\\ 
+ \sum_{k=1}^2\Vert \widetilde L_k(x_0,\widetilde x,D,\tau) 
((\partial_{\widetilde x_n}\widetilde w)e^{\tau\varphi})\Vert^2_{L^2(G)}\nonumber\\
 \le C_{10}\biggl(
\int_{\partial G\setminus\widetilde {S}_0} 
\sum_{\ell=0}^1( \tau \lambda \varphi \vert\widetilde \nabla\partial^\ell
_{\widetilde x_n}\tilde w\vert^2+\tau^3 \lambda^3
 \varphi^3 \vert\partial^\ell_{\widetilde x_n}\widetilde w\vert^2)
e^{2{\tau}\varphi}  d\widetilde {S}               \nonumber\\
+ \int_{G}(\varphi \tau \vert [\mu,\widetilde L]\widetilde w\vert^2
+ \vert\partial_{\widetilde x_n}[\mu,\widetilde L]\widetilde w\vert^2 )
e^{2{\tau}\varphi} d\widetilde xdx_0\biggr).
\end{eqnarray}
Next we estimate the surface  integral on the right-hand side of (\ref{KZP1}). 
We have
$$
\int_{\partial G\setminus\widetilde {S}_0} 
\sum_{\ell=0}^1( \tau \lambda \varphi \vert\widetilde \nabla\partial^\ell
_{\widetilde x_n}\widetilde w\vert^2
+ \tau^3 \lambda^3 \varphi^3 \vert\partial^\ell_{\widetilde x_n}\widetilde w
\vert^2) e^{2{\tau}\varphi}  d\widetilde {S}
$$
$$
= \int_{(0,b)\times O(y')\times\{\epsilon\} } 
\sum_{\ell=0}^1( \tau \lambda \varphi \vert\widetilde \nabla\partial^\ell
_{\widetilde x_n}\widetilde w\vert^2
+ \tau^3 \lambda^3 \varphi^3 \vert\partial^\ell_{\widetilde x_n}\widetilde w
\vert^2) e^{2{\tau}\varphi}  d\widetilde {S}
$$
$$
+ \int_{(0,b)\times \partial O(y')\times (0,\epsilon) } 
\sum_{\ell=0}^1( \tau \lambda \varphi \vert\widetilde \nabla\partial^\ell
_{\widetilde x_n}\widetilde w\vert^2
+ \tau^3 \lambda^3 \varphi^3 \vert\partial^\ell_{\widetilde x_n}\widetilde w
\vert^2)e^{2{\tau}\varphi}  d\widetilde {S} =: \mathcal I_1+\mathcal I_2.
$$
By (\ref{Kkhmer}), there exists a constant $\delta_5>0$ such that
\begin{equation}
\mathcal I_1\le C_{11} \int_{(0,b)\times O(y')\times\{\epsilon\} }
\sum_{\ell=0}^1( \tau \lambda \varphi \vert\widetilde \nabla\partial^\ell
_{\widetilde x_n}\widetilde w\vert^2
+ \tau^3 \lambda^3 \varphi^3 \vert\partial^\ell_{\widetilde x_n}\widetilde w
\vert^2)e^{2{\tau}(\varphi_I-\delta_5)}  d\widetilde {S}.
\end{equation}
Moreover, (\ref{KPuk}) yields 
\begin{equation}
\mathcal I_2\le C_{12}\int_{(0,T)\times \partial O(y')\times (0,\epsilon) } 
\sum_{\ell=0}^1( \tau \lambda \varphi \vert\widetilde \nabla\partial^\ell
_{\widetilde x_n}\widetilde w\vert^2
+ \tau^3 \lambda^3 \varphi^3 \vert\partial^\ell_{\widetilde x_n}\widetilde w
\vert^2)e^{2{\tau}e^{-\widetilde N(x_0-a)^2+\lambda K}}  d\widetilde x'dx_0.
\end{equation}

From (\ref{KZP1}), there exists a constant $\delta_6>0$ such that
\begin{eqnarray}\label{KzZP1}\Vert f e^{\tau\varphi_I}\Vert^2
_{L^2(0,b;L^2(O(y')\times\{0\}))}
+ \Vert f e^{\tau\varphi}\Vert^2_{L^2(0,b;L^2(O(y')\times\{0\}))}   \nonumber\\
\le C_{13}(e^{2{\tau}(\varphi_I-\delta_6)} +e^{2{\tau}e^{\lambda K}}
+ \tau e^{2{\tau}\mbox{max}\{\Vert \varphi\Vert_{C^0([b-\delta,b]\times 
\ooo\Omega)}, \Vert \varphi\Vert_{C^0([0,\delta]\times \ooo\Omega)}\}}  ).
\end{eqnarray}
Then for all sufficiently large $\tau$, we obtain
\begin{eqnarray}\label{KzZP1}
\Vert f e^{\tau\varphi}\Vert^2_{L^2(0,b;L^2(O(y')\times\{0\}))}\le C_{14}
 e^{2{\tau}e^{\lambda K}}.
\end{eqnarray}
Inequality (\ref{soska}) implies 
\begin{equation}\label{sup}
e^{2\tau (e^{\lambda K}+\delta_2)}\Vert f \Vert^2_{L^2(B(\delta_3, y'))}
\le \Vert f e^{\tau\varphi}\Vert^2_{L^2(0,b;L^2(O(y')\times\{0\}))}.
\end{equation}

In (\ref{KzZP1}), letting $\tau$ go to $+\infty$  and using (\ref{sup}),
verify that $f=0$ in $(0,T)\times B(\delta_3, y')$.
The proof of the theorem is complete.
$\blacksquare$
\\

{\bf Acknowledgments. }
The second author is supported by Grant-in-Aid for Scientific Research (A)
20H00117 and Grant-in-Aid for Challenging Research (Pioneering) 21K18142, JSPS.

\end{document}